\documentclass[reqno]{amsart}
\usepackage{color}
\usepackage{amsmath}
\usepackage{amsfonts}
\usepackage{amssymb}
\usepackage{amsthm}
\usepackage{newlfont}
\usepackage{graphicx}

\usepackage[bookmarks]{hyperref} 

\usepackage{color,soul} 

\newtheorem{thm}{Theorem}[section]

\newtheorem{lem}[thm]{Lemma}
\newtheorem{prop}[thm]{Proposition}
\theoremstyle{definition}

\theoremstyle{remark}
\newtheorem{rem}[thm]{Remark}

\numberwithin{equation}{section}
\numberwithin{thm}{section}
\newtheorem*{rem*}{Remark}

\newcommand{\abs}[1]{\left\vert#1\right\vert}

\newcommand{\bR}{{\mathbb{R}}}

\newcommand{\ed}{\end {document}}

\newcounter{smalllist}

\title[nonlocal transport]{remarks on a nonlocal transport}
\author[D. Li]{Dong Li}
\address[D. Li]{Department of Mathematics, 
The Hong Kong University of Science \& Technology, Clear Water Bay,
Kowloon, Hong Kong}%
\email{madli@ust.hk}

\author[J. Rodrigo]{Jose Rodrigo} 
\address[J. Rodrigo]
{Mathematics Research Centre,
Zeeman Building,
University of Warwick,
Coventry CV4 7AL,
United Kingdom}
\email{J.Rodrigo@warwick.ac.uk}
\thanks{}



\begin{document}
\begin{abstract}
We consider  a one dimensional nonlocal transport equation and its natural multi-dimensional
analogues. By using a new pointwise inequality for the Hilbert transform, we give a short proof
of a nonlinear inequality first proved by C\'{o}rdoba, C\'{o}rdoba and Fontelos in 2005. 
We also prove several new weighted inequalities for the Hilbert transform and 
 various nonlinear versions. Some of these results generalize to a related family
of nonlocal models. 
\end{abstract}

\maketitle

\section{Introduction and main results}
In this work we consider the following nonlinear and nonlocal transport equation
\begin{align} \label{eqCCF05}
\begin{cases}
\theta_t + (\mathcal H\theta) \theta_x = - \kappa \Lambda^{\gamma} \theta,
\quad (t, x) \in (0,\infty) \times \bR, \\
\theta|_{t=0}= \theta_0,\\
\end{cases}
\end{align}
where  $\theta=\theta(t,x)$ is a scalar-valued function defined on
$[0,\infty)\times \mathbb R$, and $\mathcal H$ is the Hilbert transform defined via 
$$
\mathcal H \theta := \frac 1 {\pi} 
\operatorname{PV} \int_{-\infty}^{\infty} \frac {\theta(y)} {x-y} dy.
$$
The number $\kappa\ge 0$ is the viscosity coefficient which governs the strength of
the linear dissipation.  The dissipation term $\Lambda^{\gamma} \theta =(-\Delta)^{\gamma/2} \theta$ is defined by
using the Fourier transform as
\begin{align*}
\widehat{\Lambda^{\gamma} \theta} (\xi) = |\xi|^{\gamma} \widehat{\theta}(\xi),
\end{align*}
where $0<\gamma\le 2$. In other words $\Lambda^{\gamma}$ is the operator corresponding
to the Fourier symbol $|\xi|^{\gamma}$. When $0<\gamma<2$ and $\theta$ has suitable
regularity (for example $\theta \in C^{1,1}$), one has the representation 
\begin{align*}
\Lambda^{\gamma} \theta = C_{\gamma}
\operatorname{PV} \int_{-\infty}^{\infty}
\frac{ \theta(x) -\theta(y)} {|x-y|^{1+\gamma}} dy,
\end{align*}
where $C_{\gamma}$ is a positive constant depending only on $\gamma$.  It follows that
if $\theta$ attains its global maximum at $x=x_*$, then 
\begin{align*}
(\Lambda^{\gamma} \theta )(x_*) \ge 0.
\end{align*}
By using this and the transport nature of the equation, one has for any smooth
solution $\theta$ to \eqref{eqCCF05} the $L^{\infty}$-maximum principle:
\begin{align*}
\| \theta(t,\cdot)\|_{\infty} \le \| \theta_0 \|_{\infty}, \quad\forall\, t>0.
\end{align*}
For $\kappa>0$ and regarding $L^{\infty}$ as the threshold space, the cases $\gamma<1$, $\gamma=1$, $\gamma>1$ are called supercritical, critical and
subcritical respectively. When $\kappa=0$ the model \eqref{eqCCF05} becomes
the inviscid case and it is deeply connected with the usual two-dimensional surface
quasi-geostrophic equation (cf.  \cite{kiselev} and the references therein for some recent
results). Compared with the usual Burgers equation with fractal
dissipation,  the model  \eqref{eqCCF05} in some
sense represents the simplest case of a nonlinear transport equation with nonlocal velocity and a viscous fractional dissipation.  For some other related one dimensional hydrodynamic models  
having some connection with
the 2D quasi-geostrophic equation and the 3D Euler equation, we refer the reader to
\cite{BM02}, \cite{CCCF05}, \cite{CLM85},   \cite{Gre96}, \cite{Sak03}, \cite{Sch86}, \cite{Yang96}
and the references therein for additional results.

Concerning the model \eqref{eqCCF05}, in the inviscid case $\kappa=0$, C\'{o}rdoba, C\'{o}rdoba and
  Fontelos \cite{CCF05} first
proved the breakdown of classical solutions to \eqref{eqCCF05}
for a generic class of smooth initial data. When $\kappa>0$, they also obtained the global well posedness in the subcritical case. For the critical case, global well-posedness
can be proved by adapting the method of continuity as in \cite{kiselev}.  Blow up for the supercritical
case $0\le \gamma<1/2$ was established in \cite{LRadv08}. Currently the
case $\frac 12 \le \gamma<1$ is still open. 

For the inviscid case the proof of \cite{CCF05} is based on an ingenious inequality:
\begin{align}\label{ine_tmp_001}
-\int_{\mathbb R}
\frac {\mathcal H f \cdot f_x} {x^{1+\delta}} dx
\ge C_{\delta} \int_{\mathbb R}
\frac{ f(x)^2} {x^{2+\delta}} dx,
\end{align}
where $-1<\delta<1$, $C_{\delta}>0$ is a constant depending only on $\delta$, and $f$
is an even bounded smooth (not necessarily decaying) function on $\mathbb R$ with $f(0)=0$. In the blow-up proof
the inequality \eqref{ine_tmp_001} is applied to $f(x)=\theta(0)-\theta(x)$ and thus $f$ in general does not decay at the spatial infinity.
The proof of \eqref{ine_tmp_001}
in \cite{CCF05} uses Mellin transform and complex analysis. A natural question is
whether one can give a completely real variable proof of \eqref{ine_tmp_001}. In this
direction Kiselev (see \cite{K10}) showed that for any even bounded $C^1$ function $f$ with $f(0)=0$ and 
$f^{\prime}\ge 0$ for $x>0$, the following inequality (see Proposition 26 therein)
\begin{align} \label{ine_tmp_002}
- \int_0^1 
\frac {\mathcal H f (x) f^{\prime}(x) f(x)^{p-1}}
{x^{\sigma}} dx
\ge C_0 \int_0^1 \frac {f(x)^{p+1}} {x^{1+\sigma}} dx,
\end{align}
where $p\ge 1$, $\sigma>0$ and $C_0$ is a positive constant depending on $p$ and $\sigma$. Later in \cite{SV16} Silvestre and Vicol gave four elegant proofs for the inviscid case (one should note
that the definition of the Hilbert transform $\mathcal H$ used in \cite{SV16} differs from the usual
one by a minus sign! See formula (1.2) therein).  The purpose of this paper is to revisit the model
\eqref{eqCCF05} and give several new and elementary proofs which are all real variable based. 
In Section 2 we first derive a new point-wise inequality (see Proposition \ref{prop_lower})
 for the Hilbert transform acting on even and non-increasing (on $(0,\infty)$) functions on
 $\mathbb R$, and then we show the C\'{o}rdoba-C\'{o}rdoba-Fontelos inequality by a simple
application of Hardy's inequality.  We also present several simplified arguments whose byproduct lead to a simple
proof of the Kiselev inequality \eqref{ine_tmp_002} and further improvements (in particular
we disprove the Kiselev inequality without the monotonicity constraint). In Section 3 we generalize the argument to dimensions 
$n\ge 2$ which works for the generalized surface quasi-geostrophic equations considered
in \cite{LRcmp09, Dong14, DL08}. Note that the blow-up proof here covers the full range of the
generalized surface quasi-geostrophic model.  In Section 4 we give another proof which works for general  functions having even symmetry
(note necessarily monotone decaying) for the Hilbert model case. In Section 5 we generalize
the argument to more general $\alpha$-patch type models.

\subsection*{Acknowledgements.}
 D. Li was supported in part by Hong Kong RGC grant GRF 16307317
 and 16309518. J. Rodrigo was supported in part by European Research Council, ERC Consolidator Grant no. 616797.


\section{radial decreasing for dimension $n=1$}

We shall use (often without explicit mentioning) the following Hardy's inequality.
\begin{lem}[Hardy] \label{lem_hardy}
If $1\le p<\infty$, $\tilde r\ge 0$ and $f$ is a non-negative measurable function on $(0,\infty)$.
Then
\begin{align*}
\int_0^{\infty} F(x)^p x^{p-\tilde r-3} dx \le \left(\frac p {\tilde r}\right)^p 
\int_0^{\infty} f(t)^p t^{p-\tilde r-1} dt,
\end{align*}
where $F(x) =\int_0^x f(t) dt$.
\end{lem}
\begin{proof}
See pp. 35 of \cite{Ziemer89}. Note that the $F(x)$ defined therein has an extra $1/x$ factor.
\end{proof}
\begin{prop}[A lower bound for Hilbert transform] \label{prop_lower}
Let $g$: $\mathbb R \to \mathbb R$ be an even continuously differentiable function which is non-increasing
on $[0,\infty)$. Assume $g^{\prime} \in L^1 \cap L^{\infty}$.
Then for any $0<x<\infty$, 
\begin{align*}
(\mathcal H g)(x) \ge \frac 2 {\pi} \cdot \frac 1 x \int_0^x ( g(y) - g(x) )dy.
\end{align*}
\end{prop}
\begin{rem} \label{rem_add_00a00}
For $f$ even, continuously differentiable and \emph{non-decreasing} on $[0, \infty)$ with
$f^{\prime}\in L^1\cap L^{\infty}$, we have the inequality 
\begin{align*}
- (\mathcal H f)(x) \ge \frac 2 {\pi}
\cdot \frac 1 x  \int_0^x ( f(x) -f(y)) dy, \quad \forall\, 0<x<\infty.
\end{align*}
\end{rem}
\begin{proof}
Since $g$ is even and $g^{\prime}\le 0$ on $[0,\infty)$, it is not difficult to check that
\begin{align*}
(\mathcal Hg)(x) & = \frac 1 {\pi} \int_0^{\infty} 
\log\abs{\frac {x-y}{x+y} } g^{\prime}(y) dy  \ge \frac 1 {\pi} \int_0^x \log \abs{ \frac{1- \frac yx} {1+\frac yx} } g^{\prime}(y) dy \\
& \ge \frac 2 {\pi} \int_0^x  \frac y x \cdot (-g^{\prime}(y) ) dy,
\end{align*}
where in the last inequality we used 
\begin{align*}
- \log \abs{  \frac {1-\epsilon} {1+\epsilon} } \ge 2\epsilon \quad \text{for all $0\le \epsilon<1$}.
\end{align*}
Integration by parts then yields the result.
\end{proof}
\begin{rem*}
Another more direct proof (under the same assumptions) 
is as follows. First observe that for each $0<x<\infty$, 
\begin{align*}
\mathcal H g(x) = \frac {2x} {\pi} \int_0^{\infty} \frac {g(y)-g(x)} {x^2-y^2} dy.
\end{align*}
Thanks to monotonicity, the integrand $\frac {g(y)-g(x)} {x^2-y^2}\ge 0$ in either the
regime $y<x$ or the regime $y>x$.  Thus we can restrict the integral to the regime
$0<y<x$, and obtain 
\begin{align*}
\mathcal H g(x) &\ge \frac {2x} {\pi} \int_{0<y<x} \frac {g(y)-g(x)} {x^2-y^2} dy
\notag \\
&\ge \frac {2} {\pi } \cdot \frac 1 x \int_0^x (g(y)-g(x) ) dy.
\end{align*}

\end{rem*}

Proposition \ref{prop_lower} can now be used to establish the following  lemma
which is essentially Lemma 2.2 found in \cite{CCF05}. The original proof therein relies on
Mellin transform and positivity of certain Fourier multipliers. 
Our new proof below avoids this  and is
completely real-variable based.
For simplicity we shall make the same assumption 
on the function $g$ as in Proposition \ref{prop_lower}.

\begin{lem} \label{Jun16_Lem}
For any $-1<\delta<1$, 
\begin{align*}
-\int_0^\infty \frac {g^\prime(x) (\mathcal Hg)(x)} {x^{1+\delta}} dx \ge
C_{\delta} \int_0^\infty \frac {(g(x)-g(0))^2} {x^{2+\delta}}  dx,
\quad \text{where $C_{\delta}=\frac 1 {\pi} \cdot \frac{(1+\delta)^2}{3+\delta}$.}
\end{align*}
\end{lem}
\begin{proof}
By Proposition \ref{prop_lower},  we have
\begin{align*}
\operatorname{LHS} & \ge \frac 2 {\pi} 
\cdot \int_0^{\infty} \frac { \int_0^x (g(y)-g(x) ) dy \cdot (-g^{\prime}(x) )}
{x^{2+\delta}} dx \notag \\
& = \frac 2 {\pi} \int_0^{\infty}
\frac{ \int_0^x (f(x) - f(y) ) dy \cdot f^{\prime}(x) } {x^{2+\delta}},
\end{align*}
where $f(x) = g(0)-g(x)$. Notice that
\begin{align*}
\frac d {dx}
\Bigl( \int_0^x (f(x) -f(y)) dy \Bigr)= x f^\prime(x).
\end{align*}
Now using this and successive integration by parts gives
\begin{align*}
 & \int_0^{\infty} \frac{ \int_0^x (f(x)-f(y) ) dy f^{\prime}(x)} { x^{2+\delta}} dx \notag \\
 =& \int_0^{\infty}
 (-1) \cdot 
 \frac d {dx} (  \int_0^x (f(x)-f(y)) dy
 \frac 1 {x^{2+\delta}}) f(x) dx \notag \\
 = &\; - \int_0^{\infty} \frac{f^{\prime}(x) f(x) } {x^{1+\delta}} dx
 +(2+\delta) \int_0^{\infty}
 \frac{ \int_0^x (f(x)-f(y)) dy f(x)} {x^{3+\delta}} dx \notag \\
 =& (-\frac{1+\delta}2 +2+\delta) \int_0^{\infty}
 \frac{f(x)^2} {x^{2+\delta}} dx -(2+\delta) \int_0^{\infty} \frac{\int_0^x f(y) dy f(x)} {
 x^{3+\delta} } dx \notag \\
 = &\frac{3+\delta}2 \int_0^{\infty} \frac{(F^{\prime}(x))^2} {x^{2+\delta}} dx
 -\frac{3+\delta}2 (2+\delta) \int_0^{\infty} \frac{F(x)^2} {x^{4+\delta}} dx,
 \end{align*}
where $F(x) = \int_0^x f(y)dy$.  By Hardy's inequality, we have
\begin{align*}
\int_0^{\infty} \frac{F(x)^2} {x^{4+\delta}} dx \le  \left( \frac 2 {3+\delta} \right)^2
\int_0^{\infty} \frac{ (F^{\prime}(x))^2} {x^{2+\delta}} dx.
\end{align*}
The result then follows.
\end{proof}
\begin{rem} \label{rem2.4_00}
One can even give a direct (without using Hardy) proof as follows. Write
(after using Proposition \ref{prop_lower})
\begin{align*}
 &\frac 2 {\pi} \int_0^{\infty}
\frac{ \int_0^x (f(x) - f(y) ) dy \cdot f^{\prime}(x) } {x^{2+\delta}} \notag \\
 = & \frac 2 {\pi}
\iint_{x\ge y} \frac{ \frac 12 \frac d {dx} (  (f(x)-f(y) )^2) } {x^{2+\delta}} dx dy \notag \\
 =& \frac {2+\delta}{\pi} \iint_{x\ge y} \frac{(f(x)-f(y))^2} {x^{3+\delta} } dx dy.
\end{align*}
Now using the inequality $(a-b)^2 = a^2+ b^2-2ab \ge (1-\alpha)a^2+ (1-\frac 1{\alpha})b^2$
for any $\alpha>0$, we obtain
\begin{align*}
\iint_{x\ge y} \frac{ (f(x)-f(y))^2} {x^{3+\delta}} dx dy
\ge ( 1+ \frac 1 {2+\delta} - ( \alpha+ \frac 1 {\alpha} \cdot \frac 1 {2+\delta} ) )
\int_0^{\infty} \frac {f(x)^2} {x^{2+\delta} } dx.
\end{align*}
Optimizing in $\alpha$ then yields the inequality with a slightly inferior constant
\begin{align*}
C_{\delta} = \frac 1 {\pi} \cdot ( 3+\delta - 2 \sqrt{2+\delta}).
\end{align*}

\end{rem}

\begin{rem}
In the preceding remark, it is possible to obtain the sharper bound by using the following argument.
Noting that
\begin{align*}
\int_{x \ge y} \frac{f(y)^2} {x^{3+\delta}} dx dy = \frac 1 {2+\delta}
\int_0^{\infty} \frac{ f(y)^2} {y^{2+\delta}}dy,
\end{align*}
it suffices to treat the term 
\begin{align*}
\iint_{x\ge y} \frac {f(x) f(y)} { x^{3+\delta}} dx dy&= \int_0^{\infty} \frac{f(x)}{x^{3+\delta}}
\left(\int_0^x f(y) dy\right)dx \notag \\
&= \frac{3+\delta}2  \int_0^{\infty} \frac { (\int_0^x f(y) dy)^2} {x^{4+\delta}} dx.
\end{align*}
By Cauchy-Schwartz
\begin{align*}
\left(\int_0^x f(y) dy \right)^2 \le \int_0^x f(y)^2 y^{-p} dy \cdot \frac {x^{p+1}}{p+1}.
\end{align*}
Interchanging the integral of $dx$ and $dy$ then gives
\begin{align*}
\int_0^{\infty} \frac { (\int_0^x f(y) dy)^2} {x^{4+\delta}} dx 
&\le  \int_0^{\infty} \frac {f(y)^2}{y^p} \cdot 
\frac{1}{p+1}\left( \int_y^{\infty} x^{-3-\delta+p} dx\right)
dy \\
&= \frac 1 {(p+1) (2+\delta-p)} \int_0^{\infty} \frac {f(y)^2} {y^{2+\delta}} dy.
\end{align*}
Choosing $p=\frac{1+\delta} 2$ then yields the sharper constant
\begin{align*}
C_{\delta} = (2+\delta) \cdot ( 1+ \frac 1 {2+\delta}- 2 \cdot \frac {3+\delta}2 \cdot 
(\frac 2 {3+\delta})^2 ) = \frac {(1+\delta)^2}{3+\delta}.
\end{align*}
\end{rem}

\subsection{Proof of the Kiselev inequality}
We now sketch a simple proof of the Kiselev inequality \eqref{ine_tmp_002}.
We emphasize that this inequality is stated for \emph{nondecreasing} even functions on $\mathbb R$. 
 For illustration purposes
we first consider the simple case $p=1$. By using Proposition \ref{prop_lower} (see
Remark \ref{rem_add_00a00}), we have
\begin{align*}
-\int_0^1 \frac{\mathcal Hf(x) f^{\prime}(x)}
{x^{\sigma}} dx
&\ge 
\frac 2 {\pi} \int_0^1 \frac{ \int_0^x (f(x)-f(y)) dy f^{\prime}(x)} {x^{1+\sigma}} dx \notag \\
& = \frac 1 {\pi} \iint_{\substack{ 0<x<1\\ 0<y<x}} 
\frac { \frac d {dx} (( f(x) -f (y))^2) } {x^{1+\sigma}} dx dy \notag \\
&\ge \frac {1+\sigma} {\pi}
\iint_{\substack{ 0<x<1\\ 0<y<x}} 
\frac {  ( f(x) -f (y))^2} {x^{2+\sigma}} dx dy,
\end{align*}
where in the last step we have integrated by part in the $x$-variable and dropped the harmless
boundary terms. Note that we can also keep the boundary term and derive a sharper inequality
as it is nonnegative. Next we proceed similarly as in Remark
\ref{rem2.4_00} and derive (below we shall take $0<\alpha<1$ and specify its
value at the very end)
\begin{align*}
 & \iint_{0<y<x<1} \frac {(f(x)-f(y))^2} {x^{2+\sigma}} dx dy \notag \\
 \ge & \iint_{0<y<x<1}
 \frac{ (1-\alpha) f(x)^2+ (1-\frac 1 {\alpha}) f(y)^2} { x^{2+\sigma}} dx dy \notag \\
 = & (1-\alpha) \int_0^1 \frac {f(x)^2} {x^{1+\sigma}} dx
 +(1-\frac 1 {\alpha}) \int_0^1 f(y)^2 \frac 1 {1+\sigma} \cdot (\frac 1 {y^{1+\sigma}}-1) dy 
 \notag \\
 \ge & 
 (1-\alpha) \int_0^1 \frac {f(x)^2} {x^{1+\sigma}} dx
 +(1-\frac 1 {\alpha}) \int_0^1 f(y)^2 \frac 1 {1+\sigma} \cdot \frac 1 {y^{1+\sigma}}dy 
 \notag \\
 = & (1+\frac 1 {1+\sigma}) \int_0^1 \frac {f(x)^2} {x^{1+\sigma}} dx
 - (\alpha+\frac 1 {\alpha}\cdot \frac 1 {1+\sigma} ) \cdot \int_0^1 \frac {f(x)^2}
 {x^{1+\sigma}} dx.
 \end{align*}
Choosing $\alpha= (1+\sigma)^{-\frac 12}$ then yields the result. Note that in the second
inequality above, we used the fact that $0<\alpha<1$ so that the term $-1$ in the $y$-integral 
can be safely dropped.

Next we sketch the proof for $1<p<\infty$. We start with 
\begin{align*}
-\int_0^1 \frac{\mathcal Hf(x) f^{\prime}(x)}
{x^{\sigma}} f(x)^{p-1} dx
&\ge 
\frac 2 {\pi} \underbrace{\int_0^1 \frac{ \int_0^x (f(x)-f(y)) dy f^{\prime}(x)} {x^{1+\sigma}} 
f(x)^{p-1} dx}_{=:H_1}.
\end{align*}
Note that for any $0\le s\le 1$, we have the inequality
\begin{align*}
1-s \ge (1-s)^p.
\end{align*}
This in turn implies that (note that below $f(y)/f(x)
\le 1$ for $0<y<x$ since $f$ is nondecreasing!)
\begin{align*}
(f(x)-f(y)) f(x)^{p-1} \ge (f(x)-f(y))^p.
\end{align*}
Thus
\begin{align*}
H_1 &\ge \iint_{0<y<x<1}  \frac {\frac 1{p+1}\frac d {dx} ((f(x)-f(y))^{p+1} )} {x^{1+\sigma}} 
dx dy \notag \\
&\ge \frac {1+\sigma}{p+1} \iint_{0<y<x<1} \frac { (f(x)-f(y))^{p+1}} {x^{2+\sigma}} dx dy.
\end{align*}
Now note that for any $\beta>1$, one can find a constant $c_1>0$, depending only
on $p$ and $\beta$, such that 
\begin{align*}
(1-s)^{p+1} \ge c_1\cdot (1-\beta s^{p+1}),\qquad\forall\, 0\le s\le 1.
\end{align*}
This in turn implies that
\begin{align*}
(f(x)-f(y))^{p+1} \ge c_1\cdot(  f(x)^{p+1} - \beta f(y)^{p+1}).
\end{align*}
Using this inequality we then obtain 
\begin{align*}
H_1 \ge \;\operatorname{const} \cdot  (1-\frac {\beta}{1+\sigma}) \int_0^1 \frac {f(x)^{p+1}}
{x^{1+\sigma}} dx.
\end{align*}
Hence taking $1<\beta<1+\sigma$ (say $\beta=1+\frac{\sigma}2$) then finishes the proof for the case $p>1$.

\subsection{Further remarks}
We first point it out that,  under the assumption of \emph{monotonicity}, the
Kiselev inequality \eqref{ine_tmp_002} is stronger than the C\'{o}rdoba-C\'{o}rdoba-Fontelos inequality \eqref{ine_tmp_001}. Indeed fix any $C^1$ bounded even $f$ on
$\mathbb R$ with $f^{\prime}\ge 0$ on $(0,\infty)$, apply the Kiselev inequality
to $f_L(x)= f(\frac x L)$, and we get (after a change of variable)
\begin{align*}
-\int_0^L \frac{\mathcal H f f^{\prime} } {x^{\sigma}} dx \ge
C_0 \int_0^L \frac {f(x)^2} {x^{1+\sigma}} dx.
\end{align*}
Note that $C_0$ is independent of the parameter $L$. Sending $L$ to infinity and using the Lebesgue
Monotone Convergence Theorem (note that the integrand $-Hf \cdot f^{\prime}$ is non-negative!) then yields the C\'{o}rdoba-C\'{o}rdoba-Fontelos inequality for the whole regime $\sigma>0$. 
One should note that the same argument yields the inequality
\begin{align} \notag 
- \int_0^{\infty}
\frac {\mathcal H f (x) f^{\prime}(x) f(x)^{p-1}}
{x^{\sigma}} dx
\ge C_0 \int_0^{\infty} \frac {f(x)^{p+1}} {x^{1+\sigma}} dx,
\end{align}
where $p\ge 1$ and $C_0$ depends only on $p$ and $\sigma$. 

Finally we should point it out that in the Kiselev inequality, the assumption of monotonicity
cannot be dropped in general. In what follows we shall construct a counterexample which answers a question raised by Kiselev in \cite{K10} (see Remark $1$ on page 249 therein). 

\begin{prop} \label{prop_counter_01}
For any $\sigma>0$, 
there exists an even function $f \in C_c^{\infty}(\mathbb R)$ such that 
\begin{align*}
-\int_0^1 \frac { \mathcal H f (x) f^{\prime}(x)} {x^{\sigma}} dx <0.
\end{align*}
In particular we cannot have the Kiselev inequality \eqref{ine_tmp_002} for $p=1$ without
the monotonicity assumption.
\end{prop}
\begin{rem*}
Similarly one can do the case $1<p<\infty$, but we do not present the
details here.
\end{rem*}

\begin{lem} \label{lem_09Aug20_001}
For any $\sigma>0$, one can find $\phi_A \in C_c^{\infty}(\mathbb R)$,
even and supported in $\{0<|x|<1\}$, $\phi_B \in C_c^{\infty}(\mathbb R)$,
even and supported in $|x|>1$, such that
\begin{align*}
\int_0^{\infty} \frac {\mathcal H \phi_B (x) \cdot ( \phi_A)^{\prime}(x) }
{x^{\sigma}} dx >0.
\end{align*}
\end{lem}
\begin{proof}[Proof of Lemma \ref{lem_09Aug20_001}]
Clearly 
\begin{align*}
\int_0^{\infty} \frac {\mathcal H \phi_B (x) \cdot ( \phi_A)^{\prime}(x)}
{x^{\sigma}} dx = - \int_0^{\infty}
\partial_x( \frac {\mathcal H \phi_B} {x^{\sigma}} ) \phi_A(x) dx.
\end{align*}
Now
\begin{align*}
\partial_x( \frac {\mathcal H \phi_B} {x^{\sigma}} )
= \frac 1 {x^{\sigma}} ( \Lambda \phi_B - \sigma \cdot \frac 1 x \cdot \mathcal H \phi_B).
\end{align*}
Observe that (here we use $\phi_B$ is supported in $|x|>1$)
\begin{align*}
  &(\Lambda \phi_B)(1) - \sigma (\mathcal H \phi_B)(1) \notag \\
= &\;
\alpha_1 \int_{|y|>1} \frac {-\phi_B(y)} {(1-y)^2} dy
-\frac 2 {\pi} \sigma \int_{y>1} \frac {\phi_B(y)}  {1-y^2}  dy \notag \\
=& \; - \int_{y>1}
\bigl[ \frac{\alpha_1}{(1-y)^2} 
+ \frac{\alpha_1}{(1+y)^2}
+\frac 2 {\pi} \sigma \frac 1 {1-y^2} \bigr] \phi_B(y) dy,
\end{align*}  
where $\alpha_1>0$ is an absolute constant which appear in the definition of the nonlocal
operator $\Lambda$.  It is then clear that one can choose suitable $\phi_B$ such that
\begin{align*}
(\Lambda \phi_B)(1) - \sigma (\mathcal H \phi_B)(1)  <0.
\end{align*}
By continuity we can find $x_0<1$ sufficiently close to $1$, such that
\begin{align*}
(\Lambda \phi_B)(x_0) - \sigma\cdot \frac 1 {x_0} (\mathcal H \phi_B)(x_0)  <0.
\end{align*}
Choosing $\phi_A$ to be a suitable bump function localized around $x_0$ then yields the result. 
\end{proof}

With the help of Lemma \ref{lem_09Aug20_001}, we now complete the proof of
Proposition \ref{prop_counter_01}. Choose 
\begin{align*}
f(x) = \phi_A (x) + t \phi_B(x).
\end{align*}
Then clearly (note that below we use the fact that
$\phi_B$ is supported in $|x|>1$)
\begin{align*}
-\int_0^1 \frac { \mathcal H f (x) f^{\prime}(x)} {x^{\sigma}} dx 
&= c_1- t \int_0^1 \frac {\mathcal H \phi_B (x) \cdot ( \phi_A)^{\prime}(x)}
{x^{\sigma}} dx   \notag \\
&=c_1 -t \int_0^{\infty} \frac {\mathcal H \phi_B (x) \cdot ( \phi_A)^{\prime}(x)}
{x^{\sigma}} dx,
\end{align*}
where $c_1$ is independent of $t$. Choosing $t$ sufficiently large then yields the result.

\section{Radial decreasing for dimension $n\ge 2$}
In \cite{LRcmp09, Dong14, DL08} a family of  the generalized surface quasi-geostrophic equations were introduced and studied.  The simplest inviscid case takes the form:
\begin{align*}
\partial_t g + (\Lambda^{-\alpha} \nabla g \cdot \nabla g) =0,
\end{align*}
where $n\ge 2$, $0<\alpha<2$ and $\Lambda^{-\alpha}$ corresponds to the
Fourier multiplier $|\xi|^{-\alpha}$.  These models can be viewed as natural
generalizations of the one dimensional Hilbert-type models to higher dimensions.
In what follows we shall discuss the corresponding nonlinear inequalities in analogy
with the Hilbert transform case. 

\begin{prop} \label{prop3.1}
Let $n\ge 2$ and $0< \alpha <2$. Let $g:\, \mathbb R^n \to \mathbb R$ be a radial and non-increasing
Schwartz function. Then for any $x \ne 0$,
\begin{align*}
- (\Lambda^{-\alpha} \nabla g )(x) \cdot \frac x {|x|}
\ge C_{\alpha,n} \cdot
\frac 1 {r^{n-\alpha+1}} \int_0^r (-g^{\prime}(\rho)) \cdot \rho^n d \rho,
\end{align*}
where $r=|x|$ and $C_{\alpha,n}>0$ depends only on $(\alpha,n)$. 
\end{prop}
\begin{rem} \label{rem_prop3.1_0}
Note that for $f(x)=g(0)-g(x)$ radial and nondecreasing, we have 
\begin{align*}
 (\Lambda^{-\alpha} \nabla f )(x) \cdot \frac x {|x|}
\ge C_{\alpha,n} \cdot
\frac 1 {r^{n-\alpha+1}} \int_0^r (f^{\prime}(\rho)) \cdot \rho^n d \rho.
\end{align*}
\end{rem}
\begin{proof}
Since $g$ is radial we can assume WLOG that $x=re_n=r\cdot (0,\cdots, 0, 1)$. By using
the fact that $g^{\prime}(\rho)\le 0$, we have
\begin{align*}
-(\Lambda^{-\alpha} \partial_n g)(x)
& = C_{\alpha,n} \int_0^{\infty} \int_{|\omega|=1}
\frac{\omega_n}{|r e_n - \rho \omega|^{n-\alpha}}
\cdot (-g^{\prime}(\rho)) \rho^{n-1} d\sigma(\omega) d\rho \notag \\
& \gtrsim \int_0^r ( -g^{\prime}(\rho)) \cdot \rho^{n-1} \cdot r^{-(n-\alpha)} 
\cdot \frac{\rho} r d\rho,
\end{align*}
where we have used the simple inequality
\begin{align*}
\int_{|\omega|=1}
\frac{\omega_n} { |e_n - \epsilon \omega|^{n-\alpha}} d\sigma(\omega)
\gtrsim \epsilon,\quad \text{for $0<\epsilon <1$}.
\end{align*}

\end{proof}

\begin{lem} \label{lem_fullrange_00a}
Let $n\ge 2$, $0< \alpha <2$ and $-1<\delta <1$. Let $g:\, \mathbb R^n \to \mathbb R$ be a radial and nonincreasing
Schwartz function. Then
\begin{align*}
\int_{\mathbb R^n}
\frac{ \Lambda^{-\alpha} \nabla g \cdot \nabla g }
{|x|^{n+\delta}} dx
\ge C_{\alpha,\delta,n}
\int_{\mathbb R^n}
\frac{ (g(0) - g(x))^2}{ |x|^{n+2-\alpha+\delta}} dx,
\end{align*}
where $C_{\alpha,\delta,n}>0$ depends only on $(\alpha,\delta,n)$. 
\end{lem}

\begin{proof}
Denote $f(x)=g(0)-g(x)$. Note that $f$ is non-decreasing and $f(0)=0$. 
By Proposition \ref{prop3.1} and Remark 
\ref{rem_prop3.1_0}, we have
\begin{align*}
\operatorname{LHS}
& \gtrsim \int_0^{\infty} 
\frac{ \int_0^r f^{\prime}(\rho) \rho^n d\rho \cdot f^{\prime}(r)}
{ r^{n-\alpha+\delta+2}} dr \notag \\
& = - \int_0^{\infty}
\frac{f^{\prime}(r) f(r)} {r^{2-\alpha+\delta}} dr
+ (n-\alpha+\delta+2) 
\int_0^{\infty} \frac{ \int_0^r f^{\prime}(\rho) \rho^n d \rho f(r)}
{r^{n-\alpha+\delta+3}} d r \notag \\
& = - \frac{2-\alpha+\delta}2
\int_0^{\infty} \frac{f(r)^2 }{r^{3-\alpha+\delta}} dr
+ (n-\alpha+\delta+2) \int_0^{\infty} 
\frac{f(r)^2} { r^{3-\alpha+\delta}} dr \notag \\
& \quad -(n-\alpha+\delta+2)
\int_0^{\infty}
\frac{ n \int_0^r f(\rho) \rho^{n-1} d\rho f(r) r^{n-1}}
{ r^{n-1} \cdot r^{n-\alpha+\delta+3}} dr \notag \\
& = (n+\frac{2-\alpha+\delta}2)
\int_0^{\infty} 
\frac{ (f(r)r^{n-1} )^2} {r^{2n+1-\alpha+\delta}} dr \notag \\
& \quad -n(n+2-\alpha+\delta)
\cdot (n+\frac{2-\alpha+\delta}2)
\int_0^{\infty} 
\frac{ F(r)^2} {r^{2n+3-\alpha+\delta}} dr,
\end{align*}
where
\begin{align*}
F(r) = \int_0^r f(\rho) \rho^{n-1} d\rho.
\end{align*}
Now the result follows from Hardy's inequality (see Lemma \ref{lem_hardy} and
take $p=2$, $\tilde r=2n+2-\alpha+\delta$) since
\begin{align*}
1> n (n+2-\alpha+\delta) \cdot 
\left( \frac 2 {2n+2-\alpha+\delta} \right)^2.
\end{align*}
\end{proof}
With the help of Lemma \ref{lem_fullrange_00a} one can then complete
the blow-up proof for the full range of the
generalized surface quasi-geostrophic model considered in  \cite{LRcmp09, Dong14, DL08},
we omit further details.

\section{Another short proof for Hilbert}

\begin{lem} \label{lem4.1}
Let $g: \mathbb R \to \mathbb R$ be an even Schwartz function. Then
\begin{align*}
-\int_{\mathbb R} \frac{ g^{\prime}(x)  (Hg)(x) } {x} dx 
\ge \frac  1 {\pi} \int_{\mathbb R} \frac{(g(0)-g(x))^2} {x^2} dx.
\end{align*}
\end{lem}
\begin{proof}
By taking advantage of the even symmetry, we have
\begin{align*}
\operatorname{LHS}& = \frac 2 {\pi} \int_0^{\infty} g^{\prime}(x)
( \int_0^{\infty} \frac { g(x)-g(y)} { x^2-y^2} dy) dx \notag \\
& = \frac 1 {\pi} \lim_{\epsilon \to 0} \int_0^{\infty} ( \int_0^{\infty} \frac {1} {x^2-y^2 +i\epsilon} \frac d{dx}
( (g(x)-g(y))^2 ) dx ) dy \notag \\
&= \frac 1 {\pi} \int_0^{\infty} \frac{(g(0)-g(y))^2} {y^2} dy+
\frac 2 {\pi} \int_0^{\infty} \int_0^{\infty} 
\frac{(g(x)-g(y))^2} {(x^2-y^2)^2} \cdot x dx dy \notag \\
& = \frac  1 {\pi} \int_0^{\infty} \frac{(g(0)-g(y))^2} {y^2} dy + \frac 1 {\pi} 
\int_0^{\infty} \int_0^{\infty} \frac{(g(x)-g(y))^2}{(x-y)^2(x+y)} dx dy.
\end{align*}
\end{proof}
\begin{rem}
The constant $1/\pi$ is certainly not sharp since
\begin{align*}
\int_0^{\infty} \int_0^{\infty} \frac{(g(x)-g(y))^2}{(x-y)^2(x+y)} dx dy 
& \ge  2\int_{x > y} \frac{(g(x)-g(y))^2} { (x^2-y^2) (x-y) } dx dy \notag \\
& \ge 2 \int_{x> y} \frac{(g(x)-g(y))^2} {x^3} dx dy \notag \\
& \ge  (3-2\sqrt 2) \int_0^{\infty} \frac{(g(0)-g(x))^2}{x^2} dx.
\end{align*}
\end{rem}

 Lemma \ref{lem4.1} is not directly useful for establishing blow-ups since it involves a 
 non-integrable weight $1/x$.  The next lemma fixes this issue.
 
 \begin{lem} \label{lem4.3}
Let $g: \mathbb R \to \mathbb R$ be an even Schwartz function. Then
\begin{align*}
-\int_{\mathbb R} { g^{\prime}(x)  (Hg)(x) } \frac{e^{-x}} xdx 
\ge \frac  1 {2\pi} \int_{\mathbb R} \frac{(g(0)-g(x))^2} {x^2} dx-1000 \| g\|_{\infty}^2.
\end{align*}
\end{lem}
\begin{proof}
By using the same integration by parts argument as in Lemma \ref{lem4.1}, we get
\begin{align*}
\operatorname{LHS} &= \frac 1 {\pi} \int_0^{\infty}
\frac{(g(0)-g(y))^2}{y^2} dy+\frac 2 {\pi}
\int_0^{\infty} \int_0^{\infty}
\frac{(g(x)-g(y))^2} {(x^2-y^2)^2} \cdot x e^{-x} dx dy \notag \\
& \quad +\frac 1 {2\pi} \int_0^{\infty} \int_0^{\infty} 
\frac{(g(x)-g(y))^2} {x^2-y^2} (e^{-x}-e^{-y}) dx dy. \notag 
\end{align*}
Now note 
\begin{align*}
\int_0^{\infty} \int_{\frac x 2 <y<2x}
\Bigl|\frac{(g(x)-g(y))^2}{x^2-y^2}\cdot (e^{-x}-e^{-y}) \Bigr| dx dy 
\le 100 \|g\|_{\infty}^2.
\end{align*}
Also
\begin{align*}
&\int_0^{\infty} \int_{y\le \frac x2} \frac{(g(x)-g(y))^2} {|x^2-y^2|} 
\cdot |e^{-x}-e^{-y}| dx dy  \notag \\
 \le & \; 
\frac 8 3 \int_0^{\infty}
\int_{y\le \frac x 2} \frac{(g(0)-g(x))^2+(g(0)-g(y))^2}{x^2} e^{-y} dx dy \notag \\
 \le & \; \frac 83 
 \int_0^{\infty} \frac{(g(0)-g(x))^2} {x^2} (1-e^{-\frac x2}) dx 
 +\frac 43 \int_0^{\infty} \frac{(g(0)-g(y))^2} {y} e^{-y} dy \notag \\
 \le & \; \frac 12\int_0^{\infty} \frac {(g(0)-g(x))^2} {x^2}dx+300 \| g\|_{\infty}^2.
\end{align*}
The piece $y\ge 2x$ is estimated similarly. 
\end{proof}
 
 To handle the diffusion term, we need the following auxiliary lemma.
 
%
%
%
 
 \begin{lem} \label{lem4.4}
 Let $0<\gamma<1$. Let $g:\, \mathbb R \to \mathbb R$ be an even Schwartz function.
 Then
 
 \begin{align*}
 | \int_0^{\infty} \frac{\Lambda^{\gamma} g(0) -\Lambda^{\gamma} g(x)}
 {x} e^{-x} dx | \le C_{\gamma} \int_0^{\infty} \frac{|g(0)-g(x)|}{x^{1+\gamma}}
 \log(10+\frac 1 x)dx,
 \end{align*}
 where $C_{\gamma}>0$ is a constant depending only on $\gamma$. 
 \end{lem}
 \begin{proof}
 By using parity, we have
 \begin{align*}
 (\Lambda^{\gamma} g)(x) = C_{\gamma}^{(1)}
 \int_0^{\infty} \frac{2 g(x) -g (x-y) -g (x+y)} {y^{1+\gamma}} dy,
 \end{align*}
 where $C_{\gamma}^{(1)}>0$ is a constant depending only on $\gamma$. 
 Now
 \begin{align*}
  &\int_0^{\infty} \frac{\Lambda^{\gamma} g(0) -\Lambda^{\gamma} g(x)}
 {x} e^{-x} dx \notag \\
 = &C_{\gamma}^{(1)}
 \int_{0<x<\infty, 0<y<\infty}
 \frac{ 2(g(0)-g(y)) + g(x-y)+g(x+y)- 2g(x)} {xy^{1+\gamma}} e^{-x} dx dy \notag \\
 = & C_{\gamma}^{(1)}
 \int_{0<x<\infty, 0<y<\infty}
 \underbrace{\frac{ -2\tilde g(y) + \tilde g(x-y) + \tilde g(x+y) - 2 \tilde g(x) }
 { x y^{1+\gamma}} e^{-x} }_{=:H_1} dx dy,
 \end{align*}
 where  for simplicity we have denoted $\tilde g(x) := g(x)-g(0)$. 
 
 Case 1: $\frac 1 {10} \le \frac x y \le 10$. Clearly
 \begin{align*}
 \int_{\substack{\frac 1 {10} \le \frac x y \le 10\\ 0<x<\infty} } |H_1| dx dy \lesssim
 \int_0^{\infty} \frac{ |\tilde g(x) |} {x^{1+\gamma}} dx.
 \end{align*}
 
 Case 2: $y\ge 10x$. Obviously
 \begin{align*}
 \int_{\substack{ y\ge 10 x \\ 0<x<\infty}}
 \frac{|\tilde g(x)|} {x y^{1+\gamma}} e^{-x } dx dy 
 \lesssim \int_0^{\infty} \frac{ |\tilde g(x)|} {x^{1+\gamma}} dx.
 \end{align*}
 On the other hand,
 \begin{align*}
 & |\int_{\substack{y\ge 10 x \\ 0<x<\infty}}
 \frac{2\tilde g(y)- \tilde g (y-x) -\tilde g(y+x)} {xy^{1+\gamma}} e^{-x} dx dy | \notag \\
=& |\int_{0<x,y<\infty}
\frac{\tilde g(y) e^{-x}} {x}
\cdot ( \frac 2 {y^{1+\gamma}}
\cdot 1_{y\ge 10x}
-\frac 1 {(y+x)^{1+\gamma}} \cdot 1_{y\ge 9x}
-\frac 1 {(y-x)^{1+\gamma}} \cdot 1_{y\ge 11x} ) dx dy| \notag \\
\lesssim & \int_0^{\infty} \frac{|\tilde g(y)|} {y^{1+\gamma}} dy.
\end{align*}

Case 3: $x \ge 10 y$.  First
\begin{align*}
\int_{\substack{ x\ge 10y \\ 0<y<\infty}}
\frac{|\tilde g(y)|} {xy^{1+\gamma}} e^{-x} dx dy
\lesssim \int_0^{\infty}
\frac{|\tilde g(y)|} {y^{1+\gamma}}
\cdot \log(10+\frac 1 y) dy.
\end{align*}
On the other hand,
\begin{align*}
 & | \int_{\substack{x \ge 10 y\\ 0<y<\infty}}
 \frac{ 2\tilde g (x) -\tilde g(x+y) -\tilde g (x-y)} { x y^{1+\gamma}} e^{-x} dx d y| \notag \\
 = & | \int_{0<x,y<\infty}
 \frac {\tilde g(x)} {y^{1+\gamma}}
 \cdot ( \frac 2 x e^{-x}\cdot 1_{x \ge 10y} -
 \frac{e^{-(x-y)}}{x-y} \cdot 1_{x\ge 11y} - \frac{e^{-(x+y)}}{x+y} \cdot 1_{x\ge 9y}) dx dy
 | \notag \\
 \lesssim & \int_0^{\infty} \frac{|\tilde g(x)|} {x^{1+\gamma}} dx.
\end{align*}

 \end{proof}

 Lemma \ref{lem4.3} and \ref{lem4.4} can be used to establish blow up. Consider
 \begin{align*}
 \begin{cases}
 \partial_t \theta - H\theta \theta_x =-\Lambda^{\gamma} \theta, \\
 \theta(0,x)=\theta_0(x).
 \end{cases}
 \end{align*}
 
 \begin{thm}
 Let $0<\gamma <\frac 12$.  Let the initial data $\theta_0$ be an even Schwartz function.
 There exists a constant $A_{\gamma}>0$ depending only on $\gamma$
 such that if 
 \begin{align*}
 \int_0^{\infty} \frac{\theta_0(0) - \theta_0(x) } x e^{-x} dx \ge A_{\gamma}
 \cdot (\| \theta_0\|_{\infty} +1),
 \end{align*}
 then the corresponding solution blows up in finite time.
 \end{thm}
\begin{proof}
By using Lemma \ref{lem4.3} and \ref{lem4.4}, we compute
\begin{align*}
&\frac d {dt} \int_0^{\infty}
\frac{\theta(t,0)-\theta(t,x)} x e^{-x} dx  \notag \\
\ge & \frac 1 {2\pi}
\int_0^{\infty} \frac{(\theta(t,0) - \theta(t,x))^2} {x^2} dx - 1000 \| \theta\|_{\infty}^2
\notag \\
& \quad
 - C_{\gamma} \int_0^{\infty}
\frac{|\theta(t,0) -\theta(t,x)|} {x^{1+\gamma}} \log (10+\frac 1 x) dx.
\end{align*}
By Cauchy-Schwartz, it is clear that
\begin{align*}
2 ( \int_0^{\infty} \frac{\theta(t,0)-\theta(t,x)} {x} e^{-x} dx )^2
\le \int_0^{\infty} \frac{(\theta(t,0)-\theta(t,x))^2} {x^2} dx.
\end{align*}
Also by using Cauchy-Schwartz, we have
\begin{align*}
&\int_0^1 \frac{|\theta(t,0)-\theta(t,x)|} {x^{1+\gamma}} \log(10+\frac 1 x) dx \notag \\
\le & (\int_0^1 \frac{(\theta(t,0)-\theta(t,x))^2} {x^2} e^{-x} dx )^{\frac 12}
\cdot ( \int_0^1 \frac {e^x} {x^{2\gamma}} ( \log(10+\frac 1x))^2 dx )^{\frac 12 } \notag \\
\le & C_1 \cdot (\int_0^1 \frac{(\theta(t,0)-\theta(t,x))^2} {x^2} e^{-x} dx )^{\frac 12},
\end{align*}
where $C_1>0$ depends only on $\gamma$. Note that here we used the crucial assumption
$0<\gamma<\frac 12$ for the integral to converge. It then follows easily that
\begin{align*}
&\frac d {dt} \int_0^{\infty}
\frac{\theta(t,0)-\theta(t,x)} x e^{-x} dx  \notag \\
\ge & \frac 1 {2\pi}  ( \int_0^{\infty} \frac{\theta(t,0) -\theta(t,x)} {x} e^{-x} dx )^2
- C_2 \cdot ( \| \theta_0\|_{\infty}+1)^2,
\end{align*}
where $C_2>0$ depends only on $\gamma$.  Choosing $A_{\gamma}=\sqrt{2\pi C_2}$
then yields the result.
\end{proof}

\section{The $\alpha$-case}
Remarkably the computation in section 4 can also be generalized to the case with
drift term $\Lambda^{-\alpha} \partial_x \theta$. We shall employ the same weight $1/x$.

\begin{lem} \label{lem5.1aaa}
Let $0<\alpha<1$.  Let $g:\mathbb R \to \mathbb R$ be an even Schwartz function. Then
\begin{align*}
\int_0^{\infty} \frac{ \Lambda^{-\alpha} g^{\prime}(x) \cdot g^{\prime}(x)} x
dx \ge C_{\alpha} 
\cdot \int_0^{\infty} \frac {(g(0)-g(x))^2} {x^{3-\alpha}} dx,
\end{align*}
where $C_{\alpha}>0$ depends only on $\alpha$. Similarly for $1\le \alpha <2$, 
by writing $\Lambda^{-\alpha} \partial_x =- \Lambda^{-(\alpha-1)} \mathcal H$, we have
\begin{align*}
-\int_0^{\infty} \frac{ \Lambda^{-(\alpha-1)} \mathcal Hg(x) \cdot g^{\prime}(x)} x
dx \ge C_{\alpha} 
\cdot \int_0^{\infty} \frac {(g(0)-g(x))^2} {x^{3-\alpha}} dx.
\end{align*}
\end{lem}
\begin{rem}
The case $\alpha=1$ corresponds to $\Lambda^{-1} \partial_x = -\mathcal H$ which is the Hilbert
transform case which we have treated before.
\end{rem}

\begin{proof}
We first discuss the case $0<\alpha<1$. 
By using parity, we have
\begin{align*}
(\Lambda^{-\alpha} g^{\prime})(x)
& = C_{\alpha} \int_0^{\infty}
( \frac 1 {|x-y|^{1-\alpha}} - \frac 1 {|x+y|^{1-\alpha}}  )\cdot
\frac d {dy} ( g(y) -g(x) ) dy \notag \\
& = C_{\alpha} \cdot (-1) \cdot
\int_0^{\infty} 
\frac d {dy} ( 
\frac 1 {|x-y|^{1-\alpha}}  - \frac 1 {(x+y)^{1-\alpha}} )  (g(y) -g (x) ) dy \notag \\
& = C_{\alpha} \int_0^{\infty} h(x,y) ( g(y) -g (x) ) dy,
\end{align*}
where 
\begin{align*}
h(x,y) = \frac d {dx} ( \frac 1 {|x-y|^{1-\alpha}} + \frac 1 {(x+y)^{1-\alpha}}).
\end{align*}
Now we write
\begin{align*}
&2 \int_0^{\infty} \int_0^{\infty}
\frac { h(x,y) (g(y) -g (x) ) g^{\prime}(x) } x dx dy \notag \\
=& -\int_0^{\infty} \int_0^{\infty} 
\frac {h(x,y)} x \cdot \frac d {dx} ( (g(x) -g (y) )^2 ) dx dy \notag \\
= &- \int_0^{\infty} (  \frac {h(x,y)} x ( g(x) -g(y) )^2 \Bigr|_{x=0}^{\infty} ) dy
+\int_0^{\infty} \int_0^{\infty}  \frac d {dx} ( \frac {h(x,y)} x ) \cdot (g(x)-g(y))^2 dx dy. 
\end{align*}
It is easy to check that for some positive constant $C_1>0$ (below $y>0$),  
\begin{align*}
(\partial_x h)(0,y)= \frac {d^2} {dx^2} ( \frac 1 {|x-y|^{1-\alpha}} + 
\frac 1 {(x+y)^{1-\alpha} } ) \Bigr|_{x=0} = C_1\cdot y^{-(3-\alpha)}.
\end{align*}
Thus 
\begin{align*}
 & - \int_0^{\infty} (  \frac {h(x,y)} x ( g(x) -g(y) )^2 \Bigr|_{x=0}^{\infty} ) dy \notag \\
=  & \int_0^{\infty} (\partial_x h)(0,y) (g(0)-g(y))^2 dy 
\gtrsim \int_0^{\infty}  \frac { (g(0)-g(y))^2} {y^{3-\alpha}} dy.
\end{align*}
It remains for us to check that, for all $0<x,y<\infty$, $x\ne y$,
\begin{align*}
\frac d {dx} (\frac {h(x,y)} x) = \frac d {dx} ( \frac 1 x \frac d {dx}
(\frac 1 {|x-y|^{1-\alpha}} + \frac 1 {(x+y)^{1-\alpha} } ) ) \ge 0.
\end{align*}
By scaling, it suffices prove for all $0<x<\infty$, $x \ne 1$,
\begin{align*}
\frac d {dx} ( \frac 1 x \frac d {dx}
(\frac 1 {|x-1|^{1-\alpha}} + \frac 1 {(x+1)^{1-\alpha} } ) ) \ge 0.
\end{align*}
We now make a change of variable $x=\sqrt {t}$. Then we only need to prove
\begin{align*}
\frac {d^2} {dt^2} ( \frac {1} { |\sqrt t -1|^{1-\alpha} } +
\frac 1 { (\sqrt t +1)^{1-\alpha} }) \ge 0, \quad \forall\, 0<t<\infty,\, t\ne 1.
\end{align*}
For $1<t<\infty$, one can get positivity by direct differentiation. For $0<t<1$, one
can use the fact that the function
\begin{align*}
f(s) = (1-s)^{-(1-\alpha)} + (1+s)^{-(1-\alpha)}
\end{align*}
has a non-negative binomial expansion for $0<s<1$. 

We now turn to the case $1\le\alpha<2$. The case $\alpha=1$ is already treated before in
Section 4 so we assume $1<\alpha<2$. Set $\epsilon =\alpha-1 \in (0,1)$. Then it is
not difficult to check that
\begin{align*}
- (\Lambda^{-\epsilon} \mathcal H g)(x) = C_{\alpha} \int_0^{\infty}
h(x,y) (g(y) -g(x)) dy,
\end{align*}
where 
\begin{align*}
h(x,y) & = - ( \frac{|x-y|^{\epsilon}} {x-y} + \frac 1 {(x+y)^{1-\epsilon}} ) \notag \\
& = -\frac 1 {\epsilon} \frac d {dx} ( |x-y|^{\epsilon} +(x+y)^{\epsilon}).
\end{align*}
Clearly 
\begin{align*}
\partial_x h(0, y) = \operatorname{const} \cdot y^{-(2-\epsilon)}.
\end{align*}
It then suffices to check for all $0<t<\infty$, $t\ne 1$, 
\begin{align*}
- \frac {d^2} {dt^2} ( |\sqrt t-1|^{\epsilon} + |\sqrt t +1|^{\epsilon} ) \ge 0.
\end{align*}
Again for $t>1$ the inequality follows easily from direct differentiation. For $0<t<1$,
one just observe that for $0<s<1$, the binomial coefficients in the expansion of
\begin{align*}
 f(s) = (1+s)^{\epsilon} + (1-s)^{\epsilon} = C_0+ \sum_{k\ge 0} C_k s^{2k}
 \end{align*}
 satisfies $C_k<0$ for all $k\ge 1$. 
\end{proof}

\begin{lem} \label{lem5.3aaa}
Let $g: \mathbb R \to \mathbb R$ be an even Schwartz function. 
If $0<\alpha<1$. 
Then
\begin{align*}
\int_0^{\infty} \frac{ \Lambda^{-\alpha} g^{\prime}(x) \cdot g^{\prime}(x)} x
e^{-x}
dx \ge C_{\alpha}^{(1)}
\cdot \int_0^{\infty} \frac {(g(0)-g(x))^2} {x^{3-\alpha}} dx-  C_{\alpha}^{(2)}\| g\|_{\infty}^2,
\end{align*}
where $C_{\alpha}^{(1)}>0$, $C_{\alpha}^{(2)}>0$ are constants depending only on $\alpha$.
Similarly for $1\le \alpha <2$,  by writing 
$\Lambda^{-\alpha} \partial_x =- \Lambda^{-(\alpha-1)} \mathcal H$, we have
\begin{align*}
-\int_0^{\infty} \frac{ \Lambda^{-(\alpha-1)} \mathcal Hg(x) \cdot g^{\prime}(x)} x e^{-x}
dx \ge C_{\alpha}^{(3)} 
\cdot \int_0^{\infty} \frac {(g(0)-g(x))^2} {x^{3-\alpha}} dx
-C_{\alpha}^{(4)}\| g\|_{\infty}^2,
\end{align*}
where $C_{\alpha}^{(3)}>0$, $C_{\alpha}^{(4)}>0$ depend only on $\alpha$. 
\end{lem}
\begin{proof}
We only need to modify the proof of Lemma \ref{lem5.1aaa}.  Consider first the case
$0<\alpha<1$. 
Recall that for $x,y>0$, $x\ne y$, 
\begin{align*}
h(x,y) &  = \frac  d{ dx} 
\Bigl( \frac 1 {|x-y|^{1-\alpha}} + \frac1 { (x+y)^{1-\alpha} } \Bigr) \notag \\
& = -(1-\alpha) \cdot
\Bigl(  |x-y|^{-(2-\alpha)} \operatorname{sgn}(x-y) 
+ (x+y)^{-(2-\alpha)} \Bigr).
\end{align*}
It is not difficult to check that (below $c_i>0$ are positive constants):
\begin{align*}
 &  c_1 \int_0^{\infty} \Lambda^{-\alpha} g^{\prime}(x) \cdot g^{\prime}(x)
 \cdot \frac 1 x e^{-x} dx \notag \\
 =&\;  - 
 \int_0^{\infty} \int_0^{\infty}
 \frac {h(x,y)} x e^{-x} \frac d {dx} ( (g(x) -g(y))^2) dx dy \notag \\
 =&\;  c_2
 \int_0^{\infty} \frac { (g(0)-g(y))^2} {y^{3-\alpha}} dy 
 + \int_0^{\infty} \int_0^{\infty}
 \frac d {dx} ( \frac 1 x h(x,y) ) e^{-x} 
 (g(x)-g(y))^2 dx dy \notag \\
 & \quad - \int_0^{\infty} \int_0^{\infty}
 h(x,y) \frac 1 x e^{-x} (g(x)-g(y))^2 dx dy.
 \end{align*}
 By the computation in Lemma \ref{lem5.1aaa}, we have $\frac d {dx} ( \frac1x
 h(x,y)) \ge 0$ for any $x,y>0$, $x\ne y$. Thus we only need to estimate the 
 third term above.  Observe that for $x>y$, we have $h(x,y)<0$.  Then
 \begin{align*}
&- \frac 1 {1-\alpha} \int_0^{\infty} \int_0^{\infty}
 h(x,y) \frac 1 x e^{-x} (g(x)-g(y))^2 dx dy \notag \\
=&\;
\int_0^{\infty} \int_0^{\infty}
( \frac {\operatorname{sgn}(x-y)} { |x-y|^{2-\alpha}}
+ \frac 1 {(x+y)^{2-\alpha}} ) \frac1x e^{-x}
( g(x) - g(y) )^2 dx dy \notag \\
\ge & \;
\iint_{\frac x 2 \le y \le 2x} 
\Bigl( \frac {\operatorname{sgn}(x-y)} { |x-y|^{2-\alpha}}
+ \frac 1 {(x+y)^{2-\alpha}}  \Bigr) \frac1x e^{-x}
( g(x) - g(y) )^2 dx dy  \notag \\
& \quad + \iint_{x<\frac 12 y}
\Bigl( \frac {\operatorname{sgn}(x-y)} { |x-y|^{2-\alpha}}
+ \frac 1 {(x+y)^{2-\alpha}}  \Bigr) \frac1x e^{-x}
( g(x) - g(y) )^2 dx dy \notag \\
\ge & \;
\iint_{\frac x 2 \le y \le 2x} 
\frac {\operatorname{sgn}(x-y)} { |x-y|^{2-\alpha}}
\cdot \frac1x e^{-x}
( g(x) - g(y) )^2 dx dy  \notag \\
& \quad + \iint_{x<\frac 12 y}
\Bigl( \frac {\operatorname{sgn}(x-y)} { |x-y|^{2-\alpha}}
+ \frac 1 {(x+y)^{2-\alpha}}  \Bigr) \frac1x e^{-x}
( g(x) - g(y) )^2 dx dy \notag \\
\ge &\; -\frac 12 \iint_{\frac x 2 \le y \le 2x} 
\frac { 1} {|x-y|^{2-\alpha}}
|\frac 1 x e^{-x} -\frac 1 y e^{-y} | (g(x)-g(y))^2 dx dy \notag \\
& \qquad - d_1 \iint_{ x <\frac 12 y}
\frac 1 {y^{3-\alpha}} e^{-x} (g(x)-g(y))^2 dx dy,
\end{align*}
where $d_1>0$ is a constant depending only on $\alpha$. Now observe that
for $x,y>0$ with $x\ne y$ and $\frac x 2 \le y \le 2x$, we have
\begin{align*}
\frac { 1} {|x-y|^{2-\alpha}}
|\frac 1 x e^{-x} -\frac 1 y e^{-y} | \lesssim
\frac 1 {|x-y|^{1-\alpha}} e^{-\frac 1 {10} |x|} (1+ x^{-2}).
\end{align*}
The desired result then follows from the following string of inequalities:
\begin{align*}
& \iint_{\frac x 2 \le y \le 2x}
|x-y|^{-(1-\alpha)} e^{-\frac 1 {10} |x|} 
(g(x)-g(y))^2 dx dy \lesssim  \|g\|_{\infty}^2;\\
& \iint_{\frac x 2 \le y \le 2x}
|x-y|^{-(1-\alpha)} |x|^{-2} e^{-\frac 1 {10} |x|} 
(g(x)-g(y))^2 dx dy \lesssim  
\; \int_0^{\infty} \frac {(g(x)-g(0))^2}{x^{2-\alpha}} dx;\\
&\iint_{ x <\frac 12 y}
\frac 1 {y^{3-\alpha}} e^{-x} (g(x)-g(y))^2 dx dy
 \lesssim  
\; \int_0^{\infty} \frac {(g(x)-g(0))^2}{x^{2-\alpha}} dx;\\
&\int_0^{\infty}\frac {(g(x)-g(0))^2}{x^{2-\alpha}} dx
\le \eta \int_0^{\infty}\frac {(g(x)-g(0))^2}{x^{3-\alpha}} dx
+C_{\eta,\alpha} \| g\|_{\infty}^2,
\end{align*}
where $\eta>0$ is any small constant, and $C_{\eta,\alpha}$ depends only
on ($\eta$, $\alpha$).

The above concludes the proof for the case $0<\alpha<1$.
The case for $1\le \alpha<2$ is similar.  In that case one only needs to work with
$h(x,y)$ given by (up to an unessential positive constant)
\begin{align*}
h(x,y) = \frac{|x-y|^{\epsilon}} {x-y} + \frac 1 {(x+y)^{1-\epsilon} },
\end{align*}
where $\epsilon=\alpha-1 \in [0,1)$.  In the symmetric region 
$\frac x2 \le y \le 2x$, one uses the inequality
\begin{align*}
\frac {|x-y|^{\epsilon}}{x-y}
\cdot \left| \frac 1x e^{-x} - \frac 1 y e^{-y} \right|
\lesssim |x-y|^{\epsilon} e^{-\frac 1 {10} x} (1+x^{-2}), \quad\forall\, x\ne y.
\end{align*}
In the region $0<x<\frac 12 y$, one can use the bound
\begin{align*}
\left |\frac {|x-y|^{\epsilon}}{x-y} + \frac 1 {(x+y)^{1-\epsilon}} \right|
\lesssim y^{-(2-\epsilon)} \cdot x.
\end{align*}
We omit further details.
\end{proof}

 Lemma  \ref{lem5.3aaa} can be used to establish blow up. For simplicity, consider
 for $0<\alpha<1$, the model
 \begin{align*}
 \begin{cases}
 \partial_t \theta +(\Lambda^{-\alpha} \partial_x \theta) \cdot \partial_x \theta=0,\\ \theta(0,x)=\theta_0(x);
 \end{cases}
 \end{align*}
 and for $1\le \alpha<2$, the model
 \begin{align*}
 \begin{cases}
 \partial_t \theta - (\Lambda^{-(\alpha-1) } \mathcal H \theta) \cdot \partial_x \theta=0,\\ \theta(0,x)=\theta_0(x);
 \end{cases}
 \end{align*}
 One should check that in both cases, the symbol of the operator for the drift
 term is given by $ i |\xi|^{-\alpha} \xi$ for all $0<\alpha<2$.  Alternatively, one 
 may write both models as a single equation
 \begin{align*}
 \partial_t \theta - \Lambda^s \mathcal H \theta \cdot \partial_x \theta =0,
 \end{align*}
 where $-1<s<1$. The drift term has the symbol $i |\xi|^s \operatorname{sgn}(\xi)$ so that 
 $s$ can be identified as $1-\alpha$.

 Concerning both models, we have the following result.
 \begin{thm} \label{thm_alpha_001}
 Let $0<\alpha<2$.  Let the initial data $\theta_0$ be an even Schwartz function.
 There exists a constant $A_{\alpha}>0$ depending only on $\alpha$
 such that if 
 \begin{align*}
 \int_0^{\infty} \frac{\theta_0(0) - \theta_0(x) } {x} e^{-x} dx \ge A_{\alpha}
 \cdot (\| \theta_0\|_{\infty} +1),
 \end{align*}
 then the corresponding solution blows up in finite time.
 \end{thm}
 \begin{rem*}
 One  can also consider the model with suitable dissipation term on the right
 hand side. For simplicity we do not state such results here which can be obtained
 by using similar estimates as in the previous section.
 \end{rem*}
\begin{proof}
This follows from Lemma \ref{lem5.3aaa}. One only needs to use the simple
inequality (with respect to the measure $e^{-x} dx$ on $(0,\infty)$) which
holds for any $0<\alpha<2$:
\begin{align*}
\int_0^{\infty}
\frac {|g(x)|} {x} e^{-x} dx
\le ( \int_0^{\infty}
\frac{ g(x)^2} {x^{3-\alpha}} e^{-x} dx)^{\frac 12}
\cdot ( \int_0^{\infty} \frac 1 {x^{\alpha-1}} e^{-x} dx )^{\frac 12}.
\end{align*}
\end{proof}

\begin{rem}
Strictly speaking,  the proof of Theorem \ref{thm_alpha_001} assumed the local
wellposedness (of smooth solutions) for the generalized model.  While the focus of this work is to prove nonlinear Hilbert type inequalities (for showing finite time singularity),  for the sake of completeness we sketch the  proof of local wellposedness here in this remark. 
Consider the nontrivial case with hyper-singular velocity as follows:
\begin{align*}
\partial_t \theta - (\Lambda^s \mathcal H \theta )\partial_x \theta =0,
\end{align*}
where $0<s<1$ (the case $-1<s\le 0$ is easier). First we present formal energy estimates.
For the basic $L^2$ estimate, we have
\begin{align*}
\frac 12 \frac d {dt} ( \| \theta \|_2^2) 
\le \frac 12 \| \Lambda^{s+1} \theta \|_{\infty} \| \theta \|_2^2.
\end{align*}
Next take an  integer $m>s+\frac 32$, and compute
\begin{align}
\frac 12 \frac d {dt}
( \| \partial_x^m \theta \|_2^2)
&=\int \partial_x^m ( \Lambda^s \mathcal H \theta \partial_x \theta) \partial_x^m \theta dx\notag\\
& = \int (\partial_x^m \Lambda^s \mathcal H \theta) \partial_x \theta \partial_x^m \theta 
dx   \label{rem_last_e001a} \\
&\qquad+ \int \Lambda^s \mathcal H \theta \partial_x^{m+1} \theta \partial_x^m \theta dx 
\label{rem_last_e001b} \\
&\qquad + \text{other terms}. \notag
\end{align}
It is not difficult to check that (one may take $m\ge 3$ for simplicity, but this can be sharpened)
\begin{align*}
|\text{other terms} | \lesssim  \| \theta\|_{H^m}^3.
\end{align*}
For \eqref{rem_last_e001b} one can do integration by parts and obtain
\begin{align*}
| \eqref{rem_last_e001b}| \lesssim \| \Lambda^{s+1} \theta\|_{\infty}
\| \partial_x^m \theta \|_2^2 \lesssim \| \theta \|_{H^m}^3.
\end{align*}
To handle \eqref{rem_last_e001a}, we can use Lemma \ref{lemforlast_00} which gives
\begin{align*}
\| \Lambda^s \mathcal H( \partial_x^m \theta \partial_x \theta)
- (\Lambda^s \mathcal H \partial_x^m \theta) \partial_x \theta
\|_2 \lesssim \| \partial_x^m \theta \|_2 \| \partial_x (1-\partial_{xx})^{\frac 34} \theta
\|_2 \lesssim \| \theta\|_{H^m}^2.
\end{align*}
Thanks to the skew-symmetry of the Hilbert transform operator, we have
\begin{align*}
\int \Lambda^s \mathcal H( \partial_x^m \theta \partial_x \theta) \partial_x^m \theta
dx = - \int \partial_x^m \theta \partial_x \theta \Lambda^s \mathcal H \partial_x^m
\theta dx.
\end{align*}
We can then rewrite the original term as a commutator  and obtain
\begin{align*}
| \eqref{rem_last_e001a} | \lesssim \| \theta\|_{H^m}^3.
\end{align*}
Thus we have completed the formal energy estimate in $H^m$. We should point it out
that by using the theory in \cite{DL19} one can obtain sharp energy estimate
in $H^r$ with $r>s+\frac 32$.  However we shall not dwell on this issue here.

Finally it is worthwhile pointing it out that in order to make the above formal energy
estimates rigorous, one needs to work with the regularized system
\begin{align*}
\partial_t \theta - J_{\epsilon} ( \Lambda^s \mathcal H J_{\epsilon} \theta \partial_x
J_{\epsilon} \theta)=0,
\end{align*}
where $J_{\epsilon}$ is the usual mollifier. We leave the interested reader to check
the details.

\end{rem}

\begin{lem} \label{lemforlast_00}
Let $0<s<1$. 
For any $f$, $g\in \mathcal S(\mathbb R)$, we have
\begin{align*}
\| \Lambda^s \mathcal H (fg) - (\Lambda^s \mathcal H f) g \|_2 \lesssim
\| f \|_2 \| |\xi|^s \hat g (\xi) \|_{L_{\xi}^1} \lesssim
\| f\|_2 \| (1-\partial_{xx})^{\frac 34} g \|_2.
\end{align*}
\end{lem}
\begin{rem}
Of course much better results are available. For example by using the commutator
estimate in \cite{DL19} (see Corollary 1.4 and the second remark on page 26 therein), one
can even show for any $1<p<\infty$, 
\begin{align*}
\| \Lambda^s \mathcal H (fg) - (\Lambda^s \mathcal H f) g \|_p \lesssim
\| f\|_p \| \Lambda^s g \|_{\operatorname{BMO} }.
\end{align*}
However for simplicity of presentation (and for the sake of completeness), we present
the non-sharp version here.
\end{rem}
\begin{proof}[Proof of Lemma \ref{lemforlast_00}]
Since we are in $L^2$ it is convenient to work purely on the Fourier side. One can 
write 
\begin{align*}
\mathcal F ( \Lambda^s \mathcal H (fg) - (\Lambda^s \mathcal H f) g ) (\xi)
& =  -i \cdot \frac 1 {2\pi} \int 
( |\xi|^s (\operatorname{sgn}(\xi) ) - |\eta|^s (\operatorname{sgn}(\eta) ) )
\widehat f(\eta) \widehat g(\xi-\eta) d\eta.
\end{align*}
It is easy to check that (since $0<s<1$)
\begin{align*}
| |\xi|^s (\operatorname{sgn}(\xi) ) - |\eta|^s (\operatorname{sgn}(\eta) )|
\lesssim | \xi -\eta|^s.
\end{align*}
The result then easily follows from Young's inequality.
\end{proof}

\end{document}